\documentclass[12pt]{amsart}
\usepackage{amsmath,amssymb,
amsbsy,amsfonts,latexsym,amsopn,
amstext,cite,amsxtra,euscript,amscd,bm,mathabx}
\usepackage{url}

\usepackage[hmargin=3cm,vmargin=3cm]{geometry}

\usepackage[colorlinks,linkcolor=blue,
anchorcolor=blue,citecolor=blue,backref=page]{hyperref}
\usepackage{color}
\usepackage{graphics,epsfig}
\usepackage{graphicx}
\usepackage{float}
\usepackage{epstopdf}
\hypersetup{breaklinks=true}

\everymath{\displaystyle}

\usepackage[norefs,nocites]{refcheck}

\newtheorem{theorem}{Theorem}[section]
\newtheorem{lemma}[theorem]{Lemma}

\newtheorem{example}[theorem]{Example}

\newcommand{\R}{\mathbb{R}}

\newcommand{\Q}{\mathbb{Q}}
\newcommand{\N}{\mathbb{N}}

\newcommand{\spt}{\text{spt\,}}
\newcommand{\diam}{\text{diam}}

\title[]{Arithmetic structure of the exceptional set of projections}


\begin{document}

\author[C. Chen]{Changhao Chen}
\address{Center for Pure Mathematics, School of Mathematical Sciences, Anhui University, Hefei 230601, China}
\email{chench@ahu.edu.cn}

\author[Z. Miao]{Zhengyan Miao}
\address{Center for Pure Mathematics, School of Mathematical Sciences, Anhui University, Hefei 230601, China}
\email{miaozhengyan321@163.com}

\begin{abstract} 
We study the arithmetic structure of the exceptional set of projections. For any bounded subset $E\subset \mathbb{R}^d$, let 
$$
\Omega=\{\xi\in \mathbb{R}: \dim_B(E+\xi E)=\dim_B E\}.
$$ 
We prove that either $\Omega=\{0\}$ or $\Omega$ is a subfield of $\mathbb{R}$.
We show that in general the statement does not hold for Hausdorff dimension and lower box dimension. Moreover, for any $s\in (0, 1]$ and a sequence $(r_k) \subset \R$, we construct a Ahlfors $s$-regular set $E\subset \R^2$ such that for any $r_k, k\in \mathbb{N}$, we have 
\[
\overline{\dim}_B \, \{x+r_k\, y: (x, y)\in E\} <s.
\]

\end{abstract}

\maketitle

\section{Introduction}

A fundamental problem in fractal geometry is to determine how the projections affect dimension \cite{Falconer2003, Mattila2015}. For $r\in \R$ define the projection $\pi_r: \R^2\rightarrow \R$ as $\pi_r(x, y)=x+ry$. A classical result is the following Marstrand projection theorem. Suppose that  $E\subset \R^2$ is a Borel set.

(1) If $\dim_H E\le 1$, then $\dim_H \pi_r (E)=\dim_H E$ for almost all $r\in \R$.

(2) If $\dim_H E>1$, then $\pi_r(E)$ has positive Lebesgue measure for almost all $r\in \R$.

Marstrand \cite{Marstrand} proved this projection theorem in the plane. 
Mattila \cite{Mattila1975}  extended the projection theorem to higher dimensions  via Kaufman’s\cite{Kaufman} potential 
theoretic methods. There is abundant literature about various extensions of Marstrand projection theorem, see  Falconer,
Fraser, and  Jin \cite{FFJ} (2015) for more details.

Let $E\subseteq \R^2$ be a Borel set and $0\le u\le \dim_H E$. The problem of exceptional projections asks for the best dimension upper bound for the  set 
\[
\{r\in \R: \dim_H \pi_r(E)<u\}.
\] 
The investigation of this problem was initiated by Kaufman \cite{Kaufman} and later by Kaufman and Mattila \cite{Mattila1975}. 

A major step in the methodology was obtained by Bourgain \cite{Bourgain2010} who introduced  the discretized sum-product technique to investigate the problem. 
Suppose $E\subset \R^2$ is a Borel set. Then Bourgain \cite{Bourgain2010} proved that 
\[
\dim_H \{r\in \R: \dim_H \pi_r(E)\le  \dim_H E \,  /2\}=0.
\]
In general the exceptional set $\{r\in \R: \dim_H \pi_r(E)\le \dim_H E/2\}$ may not be empty. Orponen \cite[Theorem 1.5]{Orponen2015} showed an Ahlfors $1$-regular set $E\subset \R^d$ such that the set 
\[
E_H(0)=\{r\in \R: \dim_H \pi_r(E)=0\}
\]
is a dense $G_\delta$ set, in particular $E_H(0)$ is uncountable.

For the dimension bound of the exceptional set, building on a series of  works  by Orponen and Shmerkin \cite{OS2023D, OS2023}, and applying tools from harmonic analysis, Ren and Wang \cite{RenWang2023} recently obtained the following sharp bound: for any Borel set $E\subset \R^2$ and $u\le \min\{\dim_H E, 1\}$,
\[
\dim_H \{r\in \R: \dim_H \pi_r(E)\le \dim_H E/2\}\le \max\{2u-\dim_H E, 0\}.\]

Besides the size, it is of interest to determine also the geometry of the possible exceptional sets.
Mattila \cite[p86]{Mattila2015} asked the following general question about the structure of exceptional sets of projections. \emph{``What more in addition to dimension estimates could be said about the exceptional sets? .... Can something be said about their structure, for example, could smooth sets or simple self-similar sets appear as exceptional sets...}"

We intend to study the structure of the exceptional set from the arithmetic viewpoint. Let  $M_{d\times d}(\R)$ denote  the $d\times d$ real matrix ring under general addition and multiplication. For $E, F \subset \R^d$, the sum sets  $E+F$ is defined as 
\[
E+F=\{x+y: x\in E, y\in F\}.
\]
For $E\subset \R^d$ and $a\in M_{d\times d}(\R)$, denote $aE=\{ax: x\in E\}$.

\begin{theorem}
\label{thm:field}
Suppose that $E\subset\R^d$ and its box dimension exists.  Then the set
\[
\{a\in M_{d\times d}(\R): \overline{\dim}_B (E+aE)=\dim_B E\}
\]
is a subring of the matrix ring $M_{d\times d}(\R)$. In particular, let
$$
\Omega=\{\xi\in \R: \overline{\dim}_B(E+\xi E)=\dim_B E\}.
$$ 
Then either $\Omega=\{0\}$ or $\Omega$ is a subfield of $\R$ under the normal addition and multiplication. 
\end{theorem}

The following Theorem \ref{thm:example} shows that the set $\Omega$ of Theorem \ref{thm:field}  can be a real subfield of $\R$. Moreover, the statement does not hold  for Hausdorff dimension and  lower box dimension.


\begin{theorem}\label{thm:example}
For any $s\in (0,d)$, there exits a compact subset $E\subseteq [0,1]^d$ such that $\dim_HE=\dim_B E=s$ and the set $E$ has the following properties: 
    \medskip
    
(i) the rational number field  $\Q \subseteq \{\xi\in \R: \dim_B(E+\xi E)=\dim_B E\} \subsetneq
 \R$;
 \medskip

(ii) the set $\{\xi\in \R:  \dim (E+\xi E)=\dim E\}$ is not a subfield of $\R$ when   $\dim$ denotes Hausdorff dimension or lower box dimension.
\end{theorem}

For the $\Omega$ in the Theorem \ref{thm:field}, it should be already known that there exists $E\subseteq \R^d$ such that the corresponding $\Omega$ contains $0$ only.  Since we could not find an explicit reference, we include it for the sake of completeness.

\begin{example}\label{ex:>s}
For any $s\in (0, d)$, there exists a subset  $E\subset [0,1]^d$ with $\dim_HE=\dim_BE=s$ such that $\dim_H (E+\xi E)>s$ holds for all $\xi \neq 0$. 
\end{example}
\begin{proof}
For $s\in (0, d)$, there exits a subset $E\subseteq \R^d$ such that $\dim_HE=\dim_BE=s$ and $E$ has positive Fourier dimension, i.e., there is a measure $\eta$ on $E$ such that 
\[
\widehat{\eta}(x)=\int_{\R^d}e^{-2\pi i x y}d\mu(y)\lesssim (1+|x|)^{-\tau/2}, \quad \forall \, x\in \R^d
\]
for some positive $\tau>0$. Here and hereafter $X\lesssim Y$ means that there exits a constant $C$ such that $|X|\le C |Y|$. Indeed, let $E$ be a typical random fractal percolation with prefixed parameters, then by Falconer \cite[Theorem 15.2 and Proposition 15.4]{Falconer2003}, we obtain $\dim_H E=\dim_B E=s$. Moreover, Suomala and Shmerkin \cite[Theorem 14.1]{ShmerkinSuomala2018}  shows that $E$ has positive Fourier dimension.

Since  $\dim_H E=s$, for any $\varepsilon>0$ by Mattila \cite[Theorem 2.8 and Theorem 3.10]{Mattila2015}, there exits a Radon measure $\mu$ on $E$ (i.e., $\spt \mu \subseteq E$) such that 
\[
I_{s-\varepsilon}(\mu)=\int_{R^d} |x|^{s-\varepsilon-d} |\widehat{\mu}(x)|^2 dx<\infty.
\]
Let $\xi\neq 0$, then $\nu=\eta\circ S_\xi^{-1}$ is a measure on $\xi E$ where $S_\xi(x)=\xi x$. Note that $\widehat{\nu}(x)\lesssim (1+|x|)^{-\tau/2}$ holds for all $\xi \in \R^d$.  It follows that 
\[
I_{s-\varepsilon+\tau}(\mu*\nu)=\int_{R^d} |x|^{s-\varepsilon+\tau-d} |\widehat{\mu}(x)|^2|\widehat{\nu}(x)|^2dx\lesssim I_{s-\varepsilon}(\mu)<\infty.
\]
Observe that $  \spt (\mu*\nu) \subseteq E+\xi E$, by \cite[Theorem 2.8]{Mattila2015}, we derive that 
\[
\dim_H(E+\xi E)\ge s+\tau-\varepsilon>s.
\]
The last inequality holds for small $\varepsilon>0$, and this finishes the proof. 
\end{proof}

For Marstrand's projection theorem, if the initial set $E\subseteq \R^2$ has additional structure, then  stronger version of Marstrand's theorem is known, see  Falconer, Fraser, and Ji \cite{FFJ} for more details. Recently,  Orponen \cite{OS2024} proved the following result for the projections of Ahlfors regular sets: 
let $E\subset \R^2$ be Ahlfors $s$-regular for some $s\in [0,2]$. Then $\dim_H E_B=0$, where 
\[
E_B=\{r\in \R: \overline{\dim}_B \pi_r(E)<\min\{s, 1\}\}.
\]
We remark that Wu \cite{Wu} has announced that $E_B$ is countable when $E$ is Ahlfors regular. The following result shows that Wu's result is sharp. 

\begin{theorem}
\label{thm:Ahlfors}
Let $s\in (0,1]$ and $(r_k)\subset \R$ be a sequence of real numbers. Then there exists a Ahlfors $s$-regular set $E\subset \R^2$ such that 
\[
\overline{\dim}_B \pi_{r_k}(E)<s, \quad \forall k\in \N.
\]   
\end{theorem}

Recall that a closed subset $E$ of $\R^d$  is  called  Ahlfors $s$-regular for some positive $s$, if there is a Radon measure $\mu$ on $E$ with $\spt \mu=E$, such that for any $x\in E$ and $r>0$, we have 
\[
r^s/C\le \mu(B(x, r))\le Cr^s
\]
for some positive constant $C$.

\section{Preliminary and proof of Theorem \ref{thm:field}}

\subsection{Preliminary}
Firstly we recall the definition of box dimension and Hausdorff dimension, see Falconer \cite{Falconer2003} for more details.

For any $\delta>0$ and $E\subseteq \R^d$, let $N_\delta(E)$ be the $\delta$-covering number of $E$, that is the smallest number of balls of diameter  $\delta$ that cover $E$. The lower and upper box dimensions
are defined respectively as

\[\underline{\dim}_B E=\liminf_{\delta\rightarrow 0}\frac{\log N_{\delta}(E)}{-\log\delta} \text{ \quad and \quad } \overline{\dim}_B E=\limsup_{\delta\rightarrow 0}\frac{\log N_{\delta}(E)}{-\log\delta}.\]
If $\underline{\dim}_B E= \overline{\dim}_B E$, then we denote this common value by $\dim_B E$ and call it the box dimension of $E$.

For $E\subseteq \R^d$, the Hausdorff dimension of $F$ is defined as 
\begin{equation*}
\dim_H E=\inf\left\{s>0: \forall \varepsilon>0, \exists \{U_i\}_{i\in \N}, U_i \subseteq \R^d, s.t. \sum_{i=1}^\infty \diam(U_i)^s< \varepsilon
\right\}
\end{equation*}

The relationship between box dimension and Hausdorff dimension is shown in the following  inequality, see \cite[Equation (3.17)]{Falconer2003},
\[
\dim_H E\le \underline{\dim}_B E \le \overline{\dim}_B E.
\]
Moreover, for $ F\subseteq \R$ let $F^d=\{(x_1, \dots, x_d): x_i\in F, 1\le i\le d\}$ be the Cartesian product set of $F$. If $\dim_H F=\dim_B F=t$ for some $t\in [0,1]$, then by \cite[Chapter 7]{Falconer2003},   
\begin{equation}\label{eq:product}
\dim_H F^d=\dim_B F^d=d\,t.
\end{equation}

We will need  tools from additive  combinatorics, in particular the following standard covering number version of Ruzsa triangle inequality and Pl\"unnecke-Ruzsa inequality, see Guth, Katz and Zahl \cite[Corollary 3.4 and Proposition 3.5]{Guth-Katz-Zahl}, and Orponen \cite[Remark 4.40]{OS2023}. For more details, see the thesis of He \cite[Chapter 2]{He}.

\begin{lemma}[Ruzsa triangle inequality]
Let $A, B, C$ be subsets of $\R^d$. Then 
\begin{equation*}
N_\delta(A\pm
B) \lesssim  \frac{N_\delta(A\pm
C)N_\delta(C\pm
B)}{N_\delta(C)}.
\end{equation*}
Here and hereafter, $X\lesssim Y$ means that $X\le cY$ for some absolutely constant $c$.
\end{lemma}

\begin{lemma}[Pl\"unnecke-Ruzsa inequality]
\label{lem:P}
Let $A, B_i, 1\le i\le n$ be subsets of $\R^d$. Suppose that $N_\delta(A+B_i)\le K_iN_\delta(A)$ for each $1\le i\le n$. Then 
\[
N_\delta(B_1+\ldots+B_n)\lesssim N_\delta(A) \prod_{i=1}^{n} K_i.
\]
\end{lemma}

For more details on additive combinatorics, see Tao and Vu \cite{Tao and Vu}, and Zhao \cite{Zhao}.

\subsection{Proof of Theorem~\ref{thm:field}}

We adapt the discrete argument of Ben Green \cite[Lemma 3.4]{Green} to our setting.

Suppose that the matrices $a, b\in \Omega$ and $a$ is invertible. It is sufficient to show (i) $a^{-1}, -a\in \Omega$;  (ii) $ab \in \Omega $; and (iii) $a+b\in \Omega$. Note that the definition of box dimension implies that 
\begin{equation}
\label{eq:basic}
N_\delta(E+xE) \leqslant  \delta^{-s+o(1)}, \quad  \delta\rightarrow 0,
\end{equation}
when $x\in \{a, b\}$. 

Since box dimension is invariant under invertible linear transform, and we obtain 
\begin{equation*}
\dim_B(a^{-1}E+E)=\dim_B(E+aE)=s.
\end{equation*}
By Ruzsa triangle inequality, we obtain  
\begin{equation}\label{eq:1}
N_\delta(E+E)\le \frac{N_\delta(E+aE) N_\delta(aE+E)}{N_\delta(aE)}\le \delta^{-s+o(1)},
\end{equation}
and 
\begin{align*}
N_\delta(E-aE)\le \frac{N_\delta(E-(-E)) N_\delta((-E)-aE)}{N_\delta(-E)}\le \delta^{-s+o(1)}, \quad \delta\rightarrow 0.
\end{align*}
Thus $\overline{\dim}_B(E-aE)\le s$. Moreover, since
\[
s= \underline{\dim}_B E \leqslant \underline{\dim}_B(E-aE)\leqslant\overline{\dim}_B(E-aE)\leqslant s,
\]
we have  $-a\in \Omega$ which proves the claim (i).

For the claim (ii),  by Ruzsa triangle inequality we obtain 
\begin{align*}
N_\delta(E+abE)&\le \frac{N_\delta(E+aE) N_\delta(aE+abE)}{N_\delta(aE)}\\
&\le \delta^{-s+o(1)}, \quad \delta\rightarrow 0.
\end{align*}
Then  $\dim_B (E+abE) =s$ which establishes (ii).

Finally, to prove the claim (iii), applying     \eqref{eq:basic}, \eqref{eq:1} and Pl\"unnecke-Ruzsa inequality (taking $A=E, B_1=E, B_2=aE, B_3=bE$),  we derive 
\begin{align*}
N_\delta(E+(a+b)E)&\le N_\delta(E+aE+bE)\\
& \le\frac{N_\delta(E+E)N_\delta(E+aE)N_\delta(E+bE)}{N_\delta(E)^{2}}\\
&\le  \delta^{-s+o(1)}, \quad \delta\rightarrow 0,
\end{align*}
and hence $\overline{\dim}_B(E+(a+b)E)\le s$. Since $\overline{\dim}_B(E+(a+b)E)\ge \dim_B E=s$, we obtain $\dim_B(E+(a+b)E)=s$ which finishes the proof the claim (iii).

\section{Proof of Theorem \ref{thm:example}}


We will need the following  partial homogeneous Cantor sets on $[0,1]$, which can be considered as an infinite generalized arithmetic progression.

{\bf Partial homogeneous Cantor sets:} Let $(m_k)\subset \N$ and $(n_k)\subset \N$ be sequence of integers with $1\le m_k<n_k$ for all $k$. Following Feng, Wen and Wu \cite{Feng}, we say that  the set 
\begin{equation}\label{eq:F}
F=F(m_k, n_k)=\left\{ \sum_{k=1}^{\infty} \frac{a_k}{n_1n_2\dots n_k} : 0\leqslant a_{k}< m_{k} \right\}
\end{equation}
is a partial homogeneous Cantor set. We remark that  Feng, Wen and Wu \cite{Feng}  used the geometric construction of the partial homogeneous which is equivalent to the above definition. Let
\[
t_1=\liminf_{k\rightarrow\infty}\frac{\log m_1m_2\dots m_k}{-\log m_{k+1}+\log n_1n_2\dots n_{k+1}}, \quad t_2= \limsup_{k\rightarrow \infty}\frac{\log m_1m_2\ldots m_k}{-\log n_1n_2\dots n_k}.
\]
Feng, Wen and Wu \cite[Lemma 2.2, Theorem 3.1]{Feng}  showed that
\begin{equation}\label{eq:dim}
\dim_H F=\underline{\dim}_B F=t_1, \quad \overline{\dim}_B F=t_2. 
\end{equation}

{\bf Notation:} In the following let  $F=F(\lceil 2^{ks/d} \rceil, 2^k)$ for some fixed $s\in (0, d)$. We intend to show that the Cartesian product set $E=F^d$
will be the desired set for Theorem \ref{thm:example}.  Firstly, by  \eqref{eq:product} and  \eqref{eq:dim}, we obtain
\begin{equation*}
 \dim_H E=\dim_B E=s.    
\end{equation*}

\begin{proof}[Proof of Theorem \ref{thm:example} (i)]
Note that any subfield of $\R$ contains rational number field $\Q$. Thus by Theorem \ref{thm:field}, it is sufficient to prove that $\dim_B (E+E)=\dim_BE$. Observe that $F+F=F(m_k', 2^k)$ where $m_k'=2\lceil 2^{ks/d} \rceil-1, k\in \N$. By \eqref{eq:dim}, we derive that 
\[
\dim_B (F+F)=\dim_H(F+F)=s/d.
\]
Observe that $E+E=(F+F)^d$, thus by \eqref{eq:product} we obtain $\dim_B(E+E)=s$ which finishes the proof.
\end{proof}

\begin{proof}[Proof of Theorem \ref{thm:example} (ii)] 
The claim (ii) of Theorem \ref{thm:example} will be a corollary of the  following Lemmas \ref{lem:C1}-\ref{lem:Eb}. In these three lemmas, we always let 
\begin{equation}\label{eq:E}
E=F^d    
\end{equation}
where $F=F(\lceil 2^{ks/d} \rceil, 2^k)$ is given by \eqref{eq:F} with $s\in (0, d)$.
\end{proof}

Recall that a subset of $\R$ is of the first Baire category if it is a countable union of nowhere dense sets; otherwise it is called of the second Baire category.  For the basic properties and various applications of Baire categories we refer to \cite{SS}.

\begin{lemma}\label{lem:C1}  Let $E$ be given by \eqref{eq:E}. Denote  $E_b=\{\xi\in \R: \underline{\dim}_B(E+\xi E)=\underline{\dim}_B E\}$, and $E_h=\{\xi\in \R: \dim_H (E+\xi E)=\dim_H E\}$.
Then both  $\R\setminus E_b$ and $\R\setminus E_h$ are of the first Baire category sets.
\end{lemma}
\begin{proof}
By \eqref{eq:dim} we have   $\dim_H E=\dim_B E=s$. Suppose that $\xi \in E_b$, then 
\[
s= \dim_H E\le  \dim_H (E+\xi E)\le \underline{\dim}_B(F+\xi F)\le s,
\]
thus $\xi\in E_h$, and hence $E_b\subseteq E_h$. Therefore, it is sufficient to prove that $\R\setminus E_b$ is of first category. 

For $m, k\in \N$, let $\varepsilon_m=1/m, \delta_k=1/2^k, m, k\in \N$, and  
\[
X_{m,k}=\left\{\xi\in \R: N_{\delta_{k}}(E+\xi E)\le \delta_k^{-s-\varepsilon_m}\right\}. 
\]
Then 
\[
E_b=\bigcap_{m=1}^\infty\bigcap_{n=1}^\infty\bigcup_{k=n}^\infty X_{m,k}.
\]

For any $\xi\in \Q$, by (i) we have  $\dim_B(E+\xi E)=s$, thus  $N_\delta(E+\xi E)\le \delta^{-s+o(1)}$, and hence  $\Q\subseteq \bigcup_{k\ge n_0} X_{m,k}$ holds for any $n_0\in \N$. Let 
\[
X_{m,k}'=\bigcup_{r\in X_{m, k}} (r-\delta_k, r+\delta_k).
\] 
Then for any $n_0$ the set $\bigcup_{k\ge n_0}X_{m, k}'$ is a dense open set. For any $m\in \N$ denote
\[
H_m=\bigcap_{n=1}^\infty\bigcup_{k=n}^\infty X_{m,k} \text{ \quad and \quad } H_m'=\bigcap_{n=1}^\infty\bigcup_{k=n}^\infty X_{m,k}'.
\]
Then $H_m\subseteq H_m'$ and $\R \setminus H_m'$ is of the first Baire category.  For any $m\ge 2$ we claim that $H_m'\subseteq H_{m-1}$. Indeed, if $\xi \in H_m'$ then there exist infinitely many $k$, such that $\xi\in X_{m, k}'$. Write $\xi=r_k+\theta_k$ where $r_k\in X_{m, k}$ and $|\theta_k|\le \delta_k$. Observe that 
\begin{align*}
N_{\delta_k}(E+\xi E)\le N_{\delta_k}(E+r_kE+\theta_k E)\le C\delta_k^{-s-\varepsilon_m},
\end{align*}
where $C$ is a constant which depends on $d$ only. It follows that   $N_{\delta_k} (E+\xi E) \le \delta_k^{-s-\varepsilon_{m-1}}$ holds for all large enough $k$. Thus $\xi\in H_m'$,
and hence 
$H_m'\subseteq H_{m-1}$. We conclude that for any $m\ge 2$ the set $\R\setminus H_{m-1}$ is of first category. Since $E_b=\bigcap_{m\in \N} H_m$, we obtain that $\R\setminus E_b$ is of fist category. 
\end{proof}

\begin{lemma}\label{lem:Eh} Let $E$ be given by \eqref{eq:E}. Denote $E_h=\{\xi\in \R: \dim_H (E+\xi E)=\dim_H E\}$. Then $E_h$ is not a subfield of $\R$. Moreover, we have  $\dim_H(E+\xi E)>s$ for some $\xi \in \R$. 
\end{lemma}
\begin{proof}
For $S\subset \R$ let $S^c=\R\setminus S$. For any $x\in \R$, by Lemma \ref{lem:C1} the set $E_h^c$ is of the first category, and thus $(x-E_h)^c$ is of the first category as well. Then by the fact that $\R$ is of the second category, we conclude that   the set 
\[
E_h\cap (x-E_h)=(E_h^c \cup (x-E_h)^c)^c\neq \emptyset,
\]
and hence $x\in (E_h+E_h)$. By the arbitrary choice of $x\in \R$, we obtain $E_h+E_h=\R$. 

Suppose to the contrary  that $E_h$ is a subfield of $\R$, then 
\begin{equation}\label{eq:all xi}
\dim_H(E+\xi E)=\dim_H E=s, \quad \forall \, \xi \in (E_h+E_h)=\R.
\end{equation}
On the other hand, observe that $E+\xi E=(F+\xi F)^d$ holds for any $\xi\in \R$. Applying \cite[Theorem 2.10]{Mattila2015}, we derive  
\begin{equation}\label{eq:E=Fd}
\dim_H (E+\xi E)\ge d \dim_H (F+\xi F).
\end{equation}
Note that $\dim_H F^2=2s/d$. Recalling that  Marstrand projection theorem claims that  for almost all $\xi \in \R$, 
\[
\dim_H (F+\xi F)=\min\{2s/d, 1\}.
\]
Combing with \eqref{eq:E=Fd} and the condition $s\in (0, d)$,  we derive that for almost  all $\xi\in\R$, 
\begin{equation}\label{eq:>s}
\dim_H (E+\xi E)\ge \min\{2s, d\}>s,
\end{equation}
which is contradict to the identity \eqref{eq:all xi}.
\end{proof}

\begin{lemma}\label{lem:Eb}
Let $E$ be given by \eqref{eq:E}. Denote  $E_b=\{\xi\in \R: \underline{\dim}_B(E+\xi E)=\underline{\dim}_B E\}$. Then $E_b$
is not a subfield of $\R$.
\end{lemma}
\begin{proof}
Applying the same argument as in the proof of Lemma \ref{lem:Eh}, and  
taking 
\[
\underline{\dim}_B (E+\xi E)\ge \dim_H (E+\xi E)>s
\]
instead of \eqref{eq:>s}, we obtain the desired result. 
\end{proof}

\section{Proof of Theorem \ref{thm:Ahlfors}}

We adapt the construction of Martin and Mattila \cite{MM} (or see Mattila \cite[Example 9.2]{Mattila1995}) for our setting.

Let $(p_1, p_2, \dots)$ be a probability vector with $p_n>0$ for all $n\in \N$. We can display the sequence  $(r_n)$ (with overlap choice) to obtain a new sequence $(\xi_n)$ such that 
\begin{equation}\label{eq:display}
\lim_{N\rightarrow \infty}\frac{1}{N}\sum_{n=1}^N {\bf 1}_{\{r_k\}}(\xi_n)=p_k, \quad \forall\, k\in \N.
\end{equation}
This can be done by applying probability method in the following way. Let $(X_n)$ be a sequence of independent identical distributed (i.i.d) random variables such that 
\[
\mathbb{P}(X_1=r_k)=p_k, \quad \forall\, k\in \N.
\]
Then for each $k\in \N$, the random variables ${\bf 1}_{\{r_k\}}(X_n), n\in \N$ are i.i.d random variables such that 
\[
\mathbb{P}({\bf 1}_{\{r_k\}}(X_n)=1)=p_k, \quad \text{and} \quad \mathbb{P}({\bf 1}_{\{r_k\}}(X_n)=0)=1-p_k.
\]
By the stronger law of large numbers \cite[Theorem 5.1.2]{Chung}, almost surely 
\begin{equation}\label{eq:law}
\lim_{N\rightarrow \infty}\frac{1}{N}\sum_{n=1}^N {\bf 1}_{\{r_k\}}(X_n)=p_k.
\end{equation}
Since there are countable $(r_k)$, we conclude that almost surely the equality \eqref{eq:law} holds for all $k\in \N$. Thus a typical sample sequence $(X_n)$ will satisfies \eqref{eq:display}. 

Now for each $\xi_k$ (direction of the projection), we intend to construct contract maps, and then iterated these maps to obtain the desired set. Let $s\in (0,1]$ and $\rho\in (0, 1/2]$ such that $2\rho^s=1$.  For $r\in \R$, let 
\[
f_{r, i}(x)=\rho x+a_{r,i}, \quad i=1, 2
\]
for some $a_{r, 1}, a_{r, 2}\in \R^2$ such that (denote  $B=B(0, 1)$)
\begin{itemize}
    \item the balls  $f_{r, 1}(B), f_{r, 2}(B) \subset B$ are interior disjoint
    \item  $\pi_r(f_{r, 1}(B))=\pi_r(f_{r, 2}(B))$. 
\end{itemize}
For any $x\in \R^2$ define $T_r(x)=\{f_{r, 1}(x), f_{r, 2}(x)\}$, and thus 
\[
T_{r}(B)=f_{r, 1}(B)\cup f_{r, 2}(B).
\]
Let 
$E_n=T_{\xi_1}\circ T_{\xi_2} \circ\ldots \circ T_{\xi_n}(B)$. Then $E_{n+1}\subset E_{n}$ for all $n\in \N$. Let $E=\bigcap_{n=1}^\infty E_n$. Clearly, $E$ is Ahlfors $s$-regular set. For any $r_k$, observe that 
\[
N_{\rho^n} (\pi_{r_k}(E))\lesssim 2^{n-n(k)},
\]
where $n(k)=\#\{1\le n\le N: \xi_n=r_k\}$. We note  that for the definition of upper (or lower) box dimension, it is sufficient to take a subsequence $\delta_n=c^n$ for some $c\in (0, 1)$. By \eqref{eq:display}, we obtain 
\[
\overline{\dim}_B \,\pi_{r_k}(E)\le (1-p_k)s,
\]
which finishes the proof.

\section*{Acknowledgement}

The authors are grateful to Tuomas Orponen, Ville Suomala and Meng Wu for helpful advice and discussions.

This work was supported by the National Natural Science Foundation of China Grant 12101002,

\end{document}